\documentclass[12pt, a4paper]{article}
\newtheorem{thm}{{\bf T}{\footnotesize \bf HEOREM}}
\newtheorem{lm}[thm]{{\bf L}{\footnotesize \bf EMMA}}
\newtheorem{cor}[thm]{{\bf C}{\footnotesize \bf OROLLARY}}
\newtheorem{pro}[thm]{{\bf P}{\footnotesize \bf ROPOSITION}}

\newtheorem{claim}{{\bf C}{\footnotesize \bf LAIM}}
\newtheorem{problem}[thm]{{\bf P}{\footnotesize \bf ROBLEM}}

\newenvironment{Proof}{\noindent {\it Proof.} }{\hbox{\rule{6pt}{6pt}}
\bigskip}
\newenvironment{Proofof}[1]
{\noindent {\it Proof of {#1}.} }{\hbox{\rule{6pt}{6pt}} \bigskip}

\usepackage[utf8]{inputenc}
\usepackage[dvipdfmx]{graphicx}
\usepackage{tikz}
\usepackage{tikzinclude}
\usepackage{amssymb}
\usepackage{float}
\usepackage{here}

\begin{document}
\sloppy
\title{Linkage problem on optimal $1$-planar graphs}

\author{
Shohei Koizumi,
\thanks{%
Graduate School of Science and Technology,
Niigata University,
8050 Ikarashi 2-no-cho, Nishi-ku, Niigata, 950-2181, Japan.
E-mail: {s-koizumi@m.sc.niigata-u.ac.jp}} 
Ryosuke Osaka, 
\thanks{%
Graduate School of Science and Technology,
Niigata University,
8050 Ikarashi 2-no-cho, Nishi-ku, Niigata, 950-2181, Japan.
E-mail: {}} 
Yusuke Suzuki
\thanks{%
Department of Mathematics, Niigata University, 
8050 Ikarashi 2-no-cho, Nishi-ku, Niigata, 950-2181, Japan.
Email: {\tt y-suzuki@math.sc.niigata-u.ac.jp}}
}
\date{}
\maketitle

\begin{abstract}
\noindent
Enami and Maezawa \cite{mnlinkedenami} give a complete characterization of $(s_1, s_2, \ldots, s_k)$-linked planar graphs for any $k$-tuple of positive integers. 
In this paper, we investigate linkage problems for optimal 1-planar graphs. 
In particular,
we show that every optimal 1-planar graph with connectivity $6$ is $(5, 5)$-linked.
Moreover, for an optimal $1$-planar graph $G$ that is not $(2,2,1)$-linked, 
we characterize disjoint vertex subsets $S_1, S_2, S_3$ in $G$
with $|S_1|=|S_2|=2$ and $|S_3|=1$ such that 
$G$ is not $\{S_1,S_2,S_3\}$-linked. 
\end{abstract}


\section{Introduction}
Throughout this paper, we consider only finite and simple graphs.
We denote the vertex set and the edge set of a graph $G$ by $V(G)$ and $E(G)$, respectively. 
For a vertex subset $S\subset V(G)$, 
we denote the induced subgraph of $G$ on $S$ by $G[S]$.
A vertex subset $S$ of a connected graph $G$ is a {\em cut} if $G-S$ is disconnected. 
For a cut $S$ of $G$, $S$ is a {\em $k$-cut} of $G$ if $|S|=k$. 
A cycle with length $k$ is a {\em $k$-cycle}. 
A cycle $C$ in $G$ is {\em separating} if $V(C)$ is a cut. 

A graph $G$ is {\em $k$-linked} 
if for any distinct $2k$ vertices $a_1, \ldots, a_k, b_1, \ldots, b_k$, 
$G$ has disjoint $k$ paths $P_1, P_2, \ldots, P_k$ 
such that $P_i$ joins $a_i$ and $b_i$ for all $i$. 
This notion has been extensively studied (see \cite{klinked1, Seymour1980}).
Let $S_1, S_2, \ldots, S_k$ be $k$ non-empty disjoint vertex subsets of $G$.  
A graph $G$ is $\{S_1, S_2, \ldots, S_k \}$-{\em linked}, 
if 
$G$ contains vertex-disjoint connected subgraphs $G_1, G_2, \ldots, G_k$ 
such that $S_i \subset V(G_i)$ for all $i$. 
Let $s_1, s_2, \ldots, s_k$ be $k$ positive integers. 
Then, $G$ is $(s_1, s_2, \ldots, s_k)$-{\em linked} 
if $G$ has at least $\sum_{i=1}^k s_i$ vertices and for any $k$ disjoint vertex subsets $S_1, S_2, \ldots, S_k$ of $G$ 
with $|S_i|= s_i$ for each $i \in \{1, 2, \ldots, k \}$,
$G$ is $\{S_1, S_2, \ldots, S_k \}$-linked. 
Note that 
a graph $G$ with at least $2k$ vertices is $(2, 2, \dots, 2)$-linked, 
i.e., $s_i$ = $2$ for all $1\leq i \leq k$, 
if and only if $G$ is $k$-linked. 
This notion was introduced by Chen et al. \cite{chen} and derived from the ``Graph Minor argument" related to a graph linkage problem. 
There are several results relating to linkages, minors, and subdivisions 
(see for example \cite{bollo, larman}).

A graph $G$ is {\em planar} if $G$ can be drawn on the sphere $S^2$ (or the plane) without edge crossings. 
Graphs already drawn on $S^2$ are {\em plane graphs}.
Let $G$ be a plane graph. 
Then a connected component of $S^2-G$, 
regarded as a topological space, is called a {\em face} of $G$.
That is, 
``a face'' in this paper is not necessarily homeomorphic to an (open) $2$-cell. 
In general, each boundary component of a face $f$ forms 
a closed walk of $G$. 
That is, the boundary of $f$ is the union of closed walks of $G$. 
In particular, if $f$ has the unique boundary component, 
then it is the {\em boundary closed walk\/} of $f$. 
A planar graph $G$ is $maximal$ if for any two non-adjacent vertices 
$x$ and $y$ of $G$, 
the graph obtained from $G$ by joining $x$ and $y$ by an edge is not planar. 
Every maximal planar graph with at least three vertices can be embedded 
on $S^2$ as a triangulation, that is, with each face bounded by a $3$-cycle. 
Mori \cite{mori200833linked} characterized $(3, 3)$-linked planar graphs. 
Furthermore, as an extension of this result, 
Enami and Maezawa \cite{mnlinkedenami} give a complete characterization of $(s_1, s_2, \ldots, s_k)$-linked planar graphs for any $k$-tuple of positive integers.

\begin{thm}[Enami and Maezawa \cite{mnlinkedenami}]\label{enamimnlinked}
Let $G$ be a planar graph, 
and let $s_1$ and $s_2$ be two integers with $2 \leq s_1 \leq s_2 \leq 4$.
Then $G$ is $(s_1, s_2)$-linked if and only if
$G$ is maximal and $(s_2+1)$-connected. 
\end{thm}


The result above also implies that if $k \ge 3$ and $s_i, s_j \ge 2$ 
for 
two distinct indices $i, j \in \{1, 2, 3, \ldots, k\}$, 
then there is no $(s_1,s_2, s_3, \ldots, s_k)$-linked planar graph
(for further details, see \cite{mnlinkedenami}).


A $1$-{\em plane graph} is a graph drawn on $S^2$ so that 
each edge crosses at most one other edge. 
An edge in a $1$-plane graph is {\em crossing\/} 
if it crosses another edge, 
and {\em non-crossing\/} otherwise. 
A $1$-{\em planar} graph is a graph that can be drawn on 
$S^2$ as a $1$-plane graph. 
The notion of $1$-planar graphs was first introduced by Ringel \cite{Ringel1965o1pg}. 
It is known that every $1$-planar graph $G$ with at least three vertices 
satisfies $|E(G)|\leq 4|V(G)|-8$ (see e.g., \cite{4n-8}). 
A $1$-planar or $1$-plane graph $G$ is {\em optimal} if it satisfies $|E(G)|= 4|V(G)|-8$. 
Briefly, 
an optimal $1$-plane graph will be abbreviated as an O1PG in this paper. 
(In fact, 
there is no need to distinguish between ``optimal $1$-planar'' 
and ``optimal $1$-plane'', since it was shown in \cite{Re} that 
every optimal $1$-planar graph 
admits a unique $1$-planar drawing in the unlabeled sense.) 
It was shown in \cite{relationship2015suzuki} that 
every O1PG contains a $4$-connected triangulation as a spanning 
subgraph. 
The following proposition can be obtained from this result and Theorem \ref{enamimnlinked}. 

\begin{pro}
Every O1PG is $(3, 3)$-linked.
\end{pro}



The aim of this paper is to characterize $(s_1, s_2, \ldots, s_k)$-linked O1PGs for certain 
$k$-tuples of positive integers. 
In \cite{ConnectivitySuzuki}, 
Fujisawa et al. discussed connectivity of O1PGs, 
and showed that every O1PG has connectivity either $4$ or $6$. 
Observe that every O1PG with connectivity $4$ (resp. $6$) is not 
$(4, n)$-linked (resp. $(6, n)$-linked) for any integer $n\geq 2$. 

The following is our first main result in the paper. 

\begin{thm}\label{thm:1-exh}
Every O1PG with connectivity $6$ is $(5, 5)$-linked.
\end{thm}

In \cite{NNS}, it was shown that  
every O1PG is obtained from a $3$-connected quadrangulation $H$ by adding a pair of crossing edges in each face of $H$. 
In other words, the spanning subgraph of an O1PG $G$, 
which consists of all non-crossing edges of $G$ is a 
$3$-connected quadrangulation. 
We denote such the quadrangulation $H$ in $G$ 
by $Q(G) (=H)$. 

Let $G$ be a $1$-plane graph. 
Then $G$ has a graph $H$ as a {\em subdrawing} 
if $H$ is a subgraph of $G$, 
and $H$ can be obtained from $G$ on $S^2$ 
by removing all the vertices and edges that are not in $H$. 
The following is our second main result in the paper. 

\begin{figure}[ht]
\centering
{\unitlength 0.1in%
\begin{picture}(24.5000,13.2800)(13.0500,-26.5800)%
\put(37.5500,-16.9000){\makebox(0,0)[lb]{$\in S_1$}}%
\put(37.5000,-20.6500){\makebox(0,0)[lb]{$\in S_2$}}%
\put(37.5000,-24.0500){\makebox(0,0)[lb]{$\in S_3$}}%
%
\special{sh 0.100}%
\special{ia 3615 2000 54 54 0.0000000 6.2831853}%
\special{pn 8}%
\special{ar 3615 2000 54 54 0.0000000 6.2831853}%
%
\special{pn 13}%
\special{ar 3614 1634 74 74 0.0000000 6.2831853}%
%
\special{sh 1.000}%
\special{ia 3614 1634 54 54 0.0000000 6.2831853}%
\special{pn 8}%
\special{ar 3614 1634 54 54 0.0000000 6.2831853}%
%
\special{pn 0}%
\special{sh 0}%
\special{pa 3555 2280}%
\special{pa 3680 2280}%
\special{pa 3680 2405}%
\special{pa 3555 2405}%
\special{pa 3555 2280}%
\special{ip}%
\special{pn 8}%
\special{pa 3555 2280}%
\special{pa 3680 2280}%
\special{pa 3680 2405}%
\special{pa 3555 2405}%
\special{pa 3555 2280}%
\special{pa 3680 2280}%
\special{fp}%
%
\special{pn 20}%
\special{pa 2145 1370}%
\special{pa 1833 1995}%
\special{fp}%
%
\special{pn 20}%
\special{pa 2926 1995}%
\special{pa 2145 1370}%
\special{fp}%
%
\special{pn 20}%
\special{pa 2145 2620}%
\special{pa 2926 1995}%
\special{fp}%
%
\special{pn 20}%
\special{pa 1364 1995}%
\special{pa 2145 2620}%
\special{fp}%
%
\special{pn 20}%
\special{pa 2145 1370}%
\special{pa 1364 1995}%
\special{fp}%
%
\special{pn 20}%
\special{pa 2457 1995}%
\special{pa 2145 1370}%
\special{fp}%
%
\special{pn 20}%
\special{pa 2145 2620}%
\special{pa 2457 1995}%
\special{fp}%
%
\special{pn 20}%
\special{pa 1833 1995}%
\special{pa 2145 2620}%
\special{fp}%
%
\special{pn 0}%
\special{sh 0}%
\special{pa 2080 1330}%
\special{pa 2205 1330}%
\special{pa 2205 1455}%
\special{pa 2080 1455}%
\special{pa 2080 1330}%
\special{ip}%
\special{pn 8}%
\special{pa 2080 1330}%
\special{pa 2205 1330}%
\special{pa 2205 1455}%
\special{pa 2080 1455}%
\special{pa 2080 1330}%
\special{pa 2205 1330}%
\special{fp}%
%
\special{pn 13}%
\special{ar 1379 1999 74 74 0.0000000 6.2831853}%
%
\special{sh 1.000}%
\special{ia 1379 1999 54 54 0.0000000 6.2831853}%
\special{pn 8}%
\special{ar 1379 1999 54 54 0.0000000 6.2831853}%
%
\special{pn 13}%
\special{ar 2444 1999 74 74 0.0000000 6.2831853}%
%
\special{sh 1.000}%
\special{ia 2444 1999 54 54 0.0000000 6.2831853}%
\special{pn 8}%
\special{ar 2444 1999 54 54 0.0000000 6.2831853}%
%
\special{sh 0.100}%
\special{ia 1835 1995 54 54 0.0000000 6.2831853}%
\special{pn 8}%
\special{ar 1835 1995 54 54 0.0000000 6.2831853}%
%
\special{sh 0.100}%
\special{ia 2905 1995 54 54 0.0000000 6.2831853}%
\special{pn 8}%
\special{ar 2905 1995 54 54 0.0000000 6.2831853}%
%
\special{sh 1.000}%
\special{ia 2154 2604 54 54 0.0000000 6.2831853}%
\special{pn 8}%
\special{ar 2154 2604 54 54 0.0000000 6.2831853}%
\end{picture}}%
\caption{A subdrawing of $Q(G)$ isomorphic to $K_{2,4}$ with specified vertices}
\label{K24}
\end{figure}
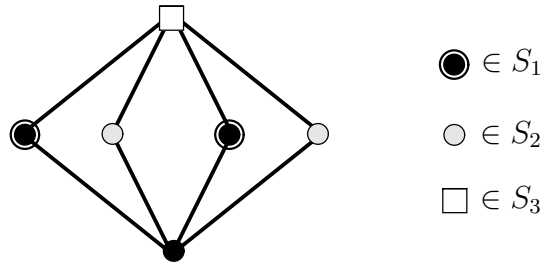

\begin{thm}\label{thm:2-exh}
Let $G$ be an O1PG,
and let $S_1$, $S_2$, and $S_3$ be disjoint vertex 
subsets of $G$ with $|S_1|=|S_2|=2$ and $|S_3|=1$. 
Then,
$G$ is $\{S_1, S_2, S_3\}$-linked
if and only if $Q(G)$ does not have a subdrawing $K$ 
such that $K$ is isomorphic to $K_{2,4}$, 
and $S_1, S_2$, and $S_3$ are arranged
on $K$ as shown in Fig.\ref{K24}. 
\end{thm}

This paper is organized as follows. 
In Section~\ref{sect:(5,5)}, 
we investigate the $(5,5)$-linkedness of O1PGs. 
We first establish some preliminary facts, 
including a general result on subtrees containing a specified set of vertices. 
We then prove the first main result of this paper. 
In Section~\ref{sect:(2,2,1)}, we consider the $(2,2,1)$-linkedness 
of O1PGs and characterize the $(2,2,1)$-linked O1PGs by determining 
the configurations of specified vertices that cannot be covered by three 
disjoint connected subgraphs.
Finally, in Section~\ref{sect:remarks}, we discuss several related topics concerning 
linkages and present some open problems.

\section{$(5,5)$-linkedness of O1PGs with connectivity 6}\label{sect:(5,5)}

In the following discussion, 
we often consider the induced subgraph of 
$Q(G)$ by a vertex subset $S$ of an O1PG $G$, 
which is $Q(G)[S]$ under our definition. 
However, when the underlying graph $G$ is clear, we use $Q[S]$ in place of $Q(G)[S]$, to simplify the notation. 
The following lemma is proved in \cite{Koizumi2025}.


\begin{lm}[Koizumi and Suzuki \cite{Koizumi2025}]\label{koizumi}
Let $G$ be an O1PG and let $S$ be a vertex subset of $G$. 
Then every face of $Q[S]$ contains at most one connected component of $G-S$. 
\end{lm}

Let $G$ be a connected graph. 
We denote the vertex subset of $G$ consisting of all vertices of 
degree $k$ by $V_k(G)$. 
For $S\subset V(G)$,  
a tree $T$ contained in $G$ as a subgraph 
is an $S$-{\em subtree} if $S\subset V(T)$ 
and all leaves of $T$ belongs to $S$. 
Note that $G$ contains at least one $S$-subtree for any $S\subset V(G)$; 
consider a spanning tree and remove vertices of degree $1$ 
that is not in $S$ recursively. 
To prove the first main theorem, we will establish the following more general statement,
which may admit other applications.

\begin{pro}\label{tree}
Let $G$ be a connected graph with $|E(G)|=|V(G)|+1$, and 
let $S$ be a subset of $V(G)$ with $|S|\geq 2$ such that 
every vertex $v$ of degree $1$ belongs to $S$. 
If every cycle in $G$ has length at least $\max\{|S|-1, 4\}$, 
and if for any two cycles $C_1$ and $C_2$,  $|V(C_1)\cap V(C_2)|\leq 2$ holds, 
then $G$ contains an $S$-subtree $T$ with $|V(T)|<|V(G)|$.  
\end{pro}

\begin{Proof} 
We use induction on $|V(G)|$. 
Let $T$ be a spanning tree of $G$, and 
denote two edges of $E(G) \setminus E(T)$ by $e_1$ and $e_2$. 
Furthermore, let 
$C_1$ (resp. $C_2$) denote the unique cycle in $T + e_1$ (resp. $T + e_2$). 
We first consider the case when $\delta(G) \geq 2$, which is actually 
the base case of the induction. 
Let $D$ be the vertex set of $G$ consisting of all vertices of degree at least $3$. 
Since $G$ has exactly $|V(G)|+1$ edges and has no vertices of degree $1$, 
we have $2|E(G)|=2(|V_2(G)| + |D| + 1) \geq 2|V_2(G)| + 3|D|$, and hence 
$|D|\leq 2$ holds. 
Therefore, we have $|S| + |D| \leq |S| + 2$. 

If $|V(C_1)\cup V(C_2)|>|S|+2$, then there is a vertex $v \in (V(C_1)\cup V(C_2))\setminus S$ 
of degree $2$. 
Thus, $G - v$ is connected, and hence $G - v$ contains a desired smaller tree. 
Thus, we assume that $|V(C_1)\cup V(C_2)| \leq |S|+2$. 
By the assumption $|C_1|, |C_2| \geq |S|-1$ and $|V(C_1)\cap V(C_2)|\leq 2$, 
we have $2|S| - 4 \leq |V(C_1)\cup V(C_2)|$, and consequently $|S| \leq 6$ holds. 
If $|S| \geq 5$, then $|V(C_1)\cup V(C_2)| \leq |S| + 2 < 2|S| - 2$ holds. 
On the other hand, if $|S| \leq 4$, then $|V(C_1)\cup V(C_2)| \leq 6$ holds by the argument above; 
note that $|C_1|, |C_2| \geq 4$ in this case by our assumption. 
That is, not depending on $|S|$, we may assume that $|V(C_1)\cap V(C_2)| \geq 1$. 
Under the condition, we have $|V(C_1)\cup V(C_2)| + 1 \leq |E(C_1)\cup E(C_2)|$. 
This implies that every component $H$ of $G - (V(C_1)\cup V(C_2))$ must be a 
tree such that $H$ is joined by exactly one edge to $C_1 \cup C_2$ by 
$|E(G)|=|V(G)|+1$. 
However, there is no such component since $\delta(G) \geq 2$. 
Thus, we now conclude that $G = C_1 \cup C_2$.

First, we assume that $|V(C_1)\cap V(C_2)| = 2$. 
Since $G$ is a two points union of two cycles, $G$ is $2$-connected . 
In addition, since $|V(G)| \geq \max \{6, 2|S|-4\}$, 
we have $|V(G)| \geq |S|+1$; 
note that if $|S|\leq 5$, then $|V(G)|\geq 6$, and 
$|V(G)| \geq 2|S|-4 \geq |S|+1$ otherwise. 
Thus, there is a vertex $v\in V(G)\setminus S$ such 
that $G - v$ is connected. 
Then, $G$ contains our desired $S$-subtree. 
Next,  assume that $|V(C_1)\cap V(C_2)| = 1$, that is $G$ is 
a one point union of $C_1$ and $C_2$. 
In this case,  $|V(G)| \geq \max \{7, 2|S|-3\}$ holds, 
and we similarly have $|V(G)| \geq |S|+2$.  
This means that $G$ has a vertex $v$ of degree $2$ 
that does not belong to $S$. 
Thus, $G - v$ contains our desired $S$-subtree.  

Then, we assume that $\delta(G) = 1$ below. 
Let $v \in V(G)$ be a vertex of degree $1$; note 
that $v \in S$. 
Let $P=v_0 \cdots v_k$ $(k \geq 1)$ denote a path with $v=v_0$ such that 
(1) $v_k$ satisfies either $v_k \in S$ or $\deg_{G}(v_k)\geq 3$, 
and (2) $v_i \notin S$ for each $\{1, \ldots, k-1\}$ (if $k\geq 2$). 
Then, put $P' = P-v_k$, and $\tilde{G} = G - P'$. 

Now, we further put $\tilde{S}=(S \setminus \{v\}) \cup \{v_k\}$. 
Then, $\tilde{G}$ and $\tilde{S}$ satisfy the conditions 
in the proposition. 
Therefore, $\tilde{G}$ contains a $\tilde{S}$-subtree $\tilde{T}$ with 
$|V(\tilde{T})|<|V(\tilde{G})|$ by the induction hypothesis. 
Then, $T=\tilde{T} \cup P$ is our desired $S$-subtree 
with $|V(T)|<|V(G)|$. 
\end{Proof}

Now, we prove the first main result in this paper. 
\bigskip

\begin{Proofof}{Theorem \ref{thm:1-exh}}
 Let $G$ be an O1PG with connectivity $6$, 
and let $S_1$, $S_2$ be disjoint subsets of $V(G)$ with 
$|S_1|=|S_2|= 5$. 
Let $G'$ be the graph obtained from $G$ by removing $S_1$. 
Since $\kappa(G) = 6$, 
$G'$ is connected. 
Thus, $G'$ contains a spanning tree. 
We choose a subtree $T$ of $G'$, which is a subdrawing of $G'$, 
so that it satisfies the following conditions $(i)$, $(ii)$, and $(iii)$. 
\begin{itemize}
\item[{$(i)$}]
$T$ contains all vertices of $S_2$, 
\item[{$(ii)$}] 
Subject to (i), the number of crossing edges of $G$ in $T$ is minimum, and 
\item[{$(iii)$}] 
Subject to (i) and (ii), $|V(T)|$ is minimum. 
\end{itemize}

\begin{claim}\label{cl.1}
All leaves of $T$ belong to $S_2$. 
\end{claim}

\begin{Proof}
Suppose that there is a leaf $v$ of $T$ which does not belong to $S_2$. 
Then, $T-v$ is a tree satisfying the condition $(i)$. 
Let $e$ be the unique edge of $T$ that is incident to $v$. 
If $e$ is a crossing edge of $G$, 
then the tree $T-v$ has fewer crossing edges than $T$, 
contradicting the condition $(ii)$. 
If $e$ is a non-crossing edge of $G$, then 
this contradicts the condition $(iii)$. 
\end{Proof}

\begin{claim}\label{cl.2}
The tree $T$ is a plane graph.
\end{claim}

\begin{Proof}
Suppose that $T$ has two edges $v_0v_2$ and $v_1v_3$ 
that are crossing each other. 
Note that since $G$ is optimal, 
there are four non-crossing edges $v_0v_1,v_1v_2, v_2v_3$, and $v_3v_0$ in $G$. 
Since $T$ is a tree,
we may assume that 
there exists the path joining $v_0$ and $v_1$ in $T$,
up to symmetry.
Then, we can obtain a new tree $T'$ from $T$ 
by removing the edge $v_1v_3$ 
and adding the edge $v_3v_0$. 
The tree $T'$ has fewer crossing edges than $T$ and 
satisfies the condition $(i)$. 
This contradicts the condition $(ii)$.
\end{Proof}

\begin{claim}\label{cl.3}
Assume that there is 
an edge $e\in E(Q[V(T)])\setminus E(T)$, 
and let $C$ be the unique cycle contained in $T + e$.
Then, every edge on $C$ is non-crossing in $G$. 
\end{claim}

\begin{Proof}
Suppose that $C$ contains a crossing edge $uv$. 
Then, $T - uv + e$ is also a tree that satisfies the condition $(i)$ (see Fig.\ref{claim3}). 
This contradicts the condition $(ii)$. 
\end{Proof}

\begin{figure}[ht]
\centering
\input{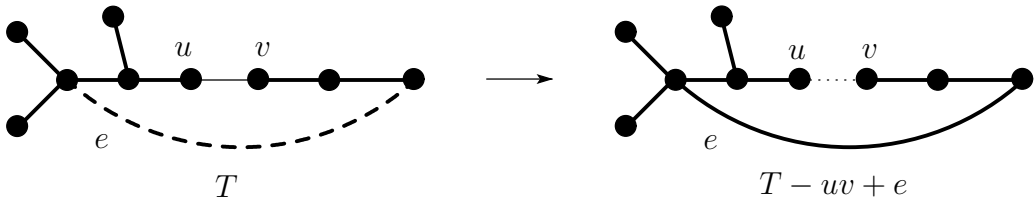}
\caption{$T$ and $T - uv + e$ in Claim $3$}
\label{claim3}
\end{figure}

\begin{claim}\label{cl.4}
Assume that there is 
an edge $e\in E(Q[V(T)])\setminus E(T)$, 
and let $C$ be the unique cycle contained in $T + e$. 
Then, $C$ has length $4$. 
Moreover, $|E(Q[V(T)]) \setminus E(T)| \le 1$, that is, 
the subgraph $Q[V(T)]$ contains at most one cycle. 
\end{claim}

\begin{Proof}
Suppose that $C$ has length at least $6$. 
Note that since $Q(G)$ is bipartite, the length of $C$ is even. 
Then,
\[
|V(C) \cap S_2|
  + |\{\, v \in V(C) : \deg_{(T+e)}(v) \geq 3 \,\}| \leq 5
\]
holds by Claim $1$.
Therefore, 
there exists a vertex $w \in V(C) \setminus S_2$ with 
$\deg_{(T+e)}(w) = 2$. 
Then, $T + e - w$ satisfies the conditions $(i)$ and $(ii)$. 
This contradicts the condition $(iii)$. 
Thus, $C$ has length $4$.

Suppose that $|E(Q[V(T)]) \setminus E(T)| \geq 2$, and let 
$e_1, e_2 \in E(Q[V(T)]) \setminus E(T)$. 
Let $C_1$ (resp. $C_2$) be the unique $4$-cycle in $T+e_1$ 
(resp. $T+e_2$). 
Then, $T\cup \{e_1, e_2 \}$ has exactly $|V(T \cup \{e_1, e_2 \})|+1$ edges. 
Suppose that $|V(C_1)\cap V(C_2)|\geq 3$. 
Since $G$ is $6$-connected and $C_1$ and $C_2$ consist only of non-crossing edges of $G$, 
each of $C_1$ and $C_2$ bounds a face of $Q(G)$. 
Thus, each of $C_1$ and $C_2$ has exactly 
two chords corresponding to crossing edges of $G$, 
contradicting the simplicity of $G$. 
Hence $|V(C_1)\cap V(C_2)|\leq 2$. 
By Proposition \ref{tree}, $T\cup \{e_1, e_2 \}$ contains a subtree $T'$ such that 
$S_2 \subset V(T')$, all leaves belong to $S_2$, and $|V(T')|<|V(T)|$. 
This contradicts the condition $(iii)$, and hence, 
$|E(Q[V(T)]) \setminus E(T)| \leq 1$. 
\end{Proof}


\begin{figure}[ht]
\centering
{\unitlength 0.1in%
\begin{picture}(47.9700,9.0100)(1.7600,-10.3700)%
%
\special{pn 20}%
\special{pa 330 790}%
\special{pa 530 590}%
\special{fp}%
%
\special{pn 20}%
\special{pa 530 590}%
\special{pa 830 690}%
\special{fp}%
%
\special{pn 8}%
\special{pa 830 690}%
\special{pa 1130 490}%
\special{fp}%
%
\special{pn 20}%
\special{pa 1130 490}%
\special{pa 1430 490}%
\special{fp}%
%
\special{pn 20}%
\special{pa 1130 490}%
\special{pa 1130 190}%
\special{fp}%
%
\special{pn 20}%
\special{pa 1430 190}%
\special{pa 1430 490}%
\special{fp}%
%
\special{pn 8}%
\special{pa 1430 490}%
\special{pa 1330 890}%
\special{fp}%
%
\special{pn 20}%
\special{pa 1730 690}%
\special{pa 1930 490}%
\special{fp}%
%
\special{pn 20}%
\special{pa 1730 690}%
\special{pa 2030 890}%
\special{fp}%
%
\special{pn 20}%
\special{pa 230 490}%
\special{pa 330 790}%
\special{fp}%
%
\special{pn 20}%
\special{pa 1330 890}%
\special{pa 1730 690}%
\special{fp}%
%
\special{pn 20}%
\special{pa 3230 790}%
\special{pa 3430 590}%
\special{fp}%
%
\special{pn 20}%
\special{pa 3430 590}%
\special{pa 3730 690}%
\special{fp}%
%
\special{pn 8}%
\special{pa 3730 690}%
\special{pa 4030 490}%
\special{dt 0.045}%
%
\special{pn 20}%
\special{pa 4030 490}%
\special{pa 4330 490}%
\special{fp}%
%
\special{pn 20}%
\special{pa 4030 490}%
\special{pa 4030 190}%
\special{fp}%
%
\special{pn 20}%
\special{pa 4030 190}%
\special{pa 4330 190}%
\special{fp}%
%
\special{pn 20}%
\special{pa 4330 190}%
\special{pa 4330 490}%
\special{fp}%
%
\special{pn 8}%
\special{pa 4330 490}%
\special{pa 4230 890}%
\special{dt 0.045}%
%
\special{pn 20}%
\special{pa 4630 690}%
\special{pa 4830 490}%
\special{fp}%
%
\special{pn 20}%
\special{pa 4630 690}%
\special{pa 4930 890}%
\special{fp}%
%
\special{pn 20}%
\special{pa 3130 490}%
\special{pa 3230 790}%
\special{fp}%
%
\special{pn 20}%
\special{pa 4230 890}%
\special{pa 4630 690}%
\special{fp}%
\put(11.3000,-11.1000){\makebox(0,0){$T$}}%
\put(40.3000,-11.1000){\makebox(0,0){$Q[V(T)]$}}%
\put(41.8000,-3.3000){\makebox(0,0){$F_1$}}%
\put(35.8000,-3.3000){\makebox(0,0){$F_2$}}%
%
\special{sh 1.000}%
\special{ia 230 490 54 54 0.0000000 6.2831853}%
\special{pn 8}%
\special{ar 230 490 54 54 0.0000000 6.2831853}%
%
\special{sh 1.000}%
\special{ia 330 775 54 54 0.0000000 6.2831853}%
\special{pn 8}%
\special{ar 330 775 54 54 0.0000000 6.2831853}%
%
\special{sh 1.000}%
\special{ia 525 600 54 54 0.0000000 6.2831853}%
\special{pn 8}%
\special{ar 525 600 54 54 0.0000000 6.2831853}%
%
\special{sh 1.000}%
\special{ia 820 680 54 54 0.0000000 6.2831853}%
\special{pn 8}%
\special{ar 820 680 54 54 0.0000000 6.2831853}%
%
\special{sh 1.000}%
\special{ia 1130 480 54 54 0.0000000 6.2831853}%
\special{pn 8}%
\special{ar 1130 480 54 54 0.0000000 6.2831853}%
%
\special{sh 1.000}%
\special{ia 1415 480 54 54 0.0000000 6.2831853}%
\special{pn 8}%
\special{ar 1415 480 54 54 0.0000000 6.2831853}%
%
\special{sh 1.000}%
\special{ia 1125 190 54 54 0.0000000 6.2831853}%
\special{pn 8}%
\special{ar 1125 190 54 54 0.0000000 6.2831853}%
%
\special{sh 1.000}%
\special{ia 1430 190 54 54 0.0000000 6.2831853}%
\special{pn 8}%
\special{ar 1430 190 54 54 0.0000000 6.2831853}%
%
\special{sh 1.000}%
\special{ia 1335 880 54 54 0.0000000 6.2831853}%
\special{pn 8}%
\special{ar 1335 880 54 54 0.0000000 6.2831853}%
%
\special{sh 1.000}%
\special{ia 1735 695 54 54 0.0000000 6.2831853}%
\special{pn 8}%
\special{ar 1735 695 54 54 0.0000000 6.2831853}%
%
\special{sh 1.000}%
\special{ia 1920 490 54 54 0.0000000 6.2831853}%
\special{pn 8}%
\special{ar 1920 490 54 54 0.0000000 6.2831853}%
%
\special{sh 1.000}%
\special{ia 2010 880 54 54 0.0000000 6.2831853}%
\special{pn 8}%
\special{ar 2010 880 54 54 0.0000000 6.2831853}%
%
\special{sh 1.000}%
\special{ia 3139 494 54 54 0.0000000 6.2831853}%
\special{pn 8}%
\special{ar 3139 494 54 54 0.0000000 6.2831853}%
%
\special{sh 1.000}%
\special{ia 3239 779 54 54 0.0000000 6.2831853}%
\special{pn 8}%
\special{ar 3239 779 54 54 0.0000000 6.2831853}%
%
\special{sh 1.000}%
\special{ia 3434 604 54 54 0.0000000 6.2831853}%
\special{pn 8}%
\special{ar 3434 604 54 54 0.0000000 6.2831853}%
%
\special{sh 1.000}%
\special{ia 3729 684 54 54 0.0000000 6.2831853}%
\special{pn 8}%
\special{ar 3729 684 54 54 0.0000000 6.2831853}%
%
\special{sh 1.000}%
\special{ia 4039 484 54 54 0.0000000 6.2831853}%
\special{pn 8}%
\special{ar 4039 484 54 54 0.0000000 6.2831853}%
%
\special{sh 1.000}%
\special{ia 4324 484 54 54 0.0000000 6.2831853}%
\special{pn 8}%
\special{ar 4324 484 54 54 0.0000000 6.2831853}%
%
\special{sh 1.000}%
\special{ia 4034 194 54 54 0.0000000 6.2831853}%
\special{pn 8}%
\special{ar 4034 194 54 54 0.0000000 6.2831853}%
%
\special{sh 1.000}%
\special{ia 4339 194 54 54 0.0000000 6.2831853}%
\special{pn 8}%
\special{ar 4339 194 54 54 0.0000000 6.2831853}%
%
\special{sh 1.000}%
\special{ia 4244 884 54 54 0.0000000 6.2831853}%
\special{pn 8}%
\special{ar 4244 884 54 54 0.0000000 6.2831853}%
%
\special{sh 1.000}%
\special{ia 4644 699 54 54 0.0000000 6.2831853}%
\special{pn 8}%
\special{ar 4644 699 54 54 0.0000000 6.2831853}%
%
\special{sh 1.000}%
\special{ia 4829 494 54 54 0.0000000 6.2831853}%
\special{pn 8}%
\special{ar 4829 494 54 54 0.0000000 6.2831853}%
%
\special{sh 1.000}%
\special{ia 4919 884 54 54 0.0000000 6.2831853}%
\special{pn 8}%
\special{ar 4919 884 54 54 0.0000000 6.2831853}%
\end{picture}}%
\vspace{5pt} 
\caption{Tree $T$ and $Q[V(T)]$}
\label{treetopological}
\end{figure}
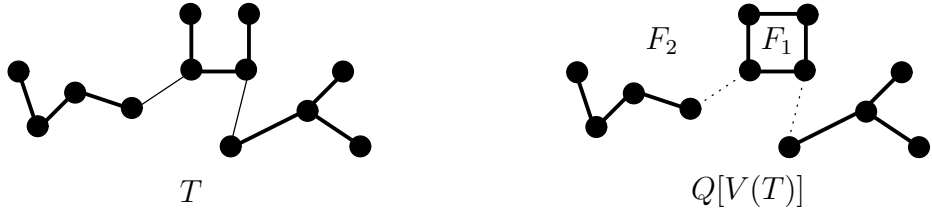

Now we consider an induced subgraph $Q[V(T)]$ of $Q(G)$.
Observe that each connected component of $Q[V(T)]$ is either acyclic or, 
if it contains exactly one cycle, then such a cycle has length exactly $4$ by 
claims~\ref{cl.3} and~\ref{cl.4}. 
%
%
%
If every connected component of $Q[V(T)]$ is acyclic, then $Q[V(T)]$ 
has a unique face. 
By Lemma \ref{koizumi},
$G - T$ is connected, and hence 
we can obtain a connected subgraph in $G - T$
containing all vertices of $S_1$. 

Therefore, we assume that there is a connected component of $Q[V(T)]$ 
that is not acyclic (see the left-hand side and the right-hand side of Fig.\ref{treetopological}). 
Then, by Claim~\ref{cl.4}, $Q[V(T)]$ has exactly two faces $F_1$ and $F_2$,
one of which, say $F_1$, is bounded by a $4$-cycle. 
Since $\kappa (G) = 6$, 
the interior of $F_1$ contains no vertices of $G$. 
Then, the interior of $F_2$ contains all vertices of $V(G) \setminus V(T)$. 
By Lemma \ref{koizumi},
$G - T$ is connected and hence 
we can obtain a connected subgraph in $G - T$
containing all vertices of $S_1$.
\end{Proofof}

\section{Characterization of $(2,2,1)$-linked O1PGs}\label{sect:(2,2,1)}


Seymour \cite{Seymour1980}, Shiloach \cite{shiloach}, and Thomassen \cite{thomassen} independently characterized $(2, 2)$-linked graphs. 
The following proposition can be obtained from the result above.  



\begin{pro}[Seymour \cite{Seymour1980}, Shiloach \cite{shiloach}, and Thomassen \cite{thomassen}]\label{seymour}
Let $G$ be a planar graph, and let $S_1=\{x_1, x_2 \}$ and 
$S_2= \{y_1, y_2 \}$ be disjoint vertex subsets of $G$. 
Then $G$ is not $\{S_1, S_2 \}$-linked if and only if 
$G$ has an embedding on $S^2$ that satisfies the followings: 

\begin{itemize}
    \item[{(i)}] There is a face $f$ whose boundary closed walk $W$ has length 
    at least $4$, and 
    \item[{(ii)}] The boundary walk $W$ contains $x_1, y_1, x_2, y_2$ in this order. 
\end{itemize}
\end{pro}

Let $G$ be a plane graph. 
Then, $G$ is a {\em near-triangulation} 
if it is $2$-connected and all faces, except at most one face, are triangles. 
A simple closed curve $\gamma$ on $S^2$ is a {\em separating $k$-curve} 
if it passes only through $k$ vertices in a plane graph $G$ 
and each of two regions separated by $\gamma$ contains at least one vertex 
of $G$.

\bigskip

\begin{Proofof}{Theorem \ref{thm:2-exh}}
Let $G$ be an O1PG, 
and let $S_1= \{r_1, r_2\}$, $S_2= \{b_1, b_2\}$, 
$S_3 = \{y\}$ be disjoint subsets of $V(G)$. 
The sufficiency is clear, and hence 
we shall prove the necessity.
Assume that $G$ is not $\{S_1, S_2, S_3\}$-linked.
Let $B$ and $W$ be the partite sets of $Q(G)$.
We may assume that $y \in W$.

Let $N_G(y) = \{u_1, u_2, \ldots, u_{2k}\}$ with $k \ge 3$,
and assume that the edges $yu_1, yu_2, \ldots, yu_{2k}$
appear in the clockwise order around $y$.
Then, we may assume that $u_i\in B$ when $i$ is even, 
and $u_i\in W$ otherwise $(i \in \{1,2, \ldots ,{2k}\})$.
Let $G'$ be the graph obtained from $G$ by removing $y$ 
(see the left-hand side of Fig.~\ref{G'andH}, 
where vertices in $B$ (resp. $W$) are shown in black (resp. white)).
Then $G'$ contains a spanning subgraph that is a near-triangulation 
whose unique non-triangular face is bounded by $C'=u_2u_4 \ldots u_{2k}$.  
Since $C'$ consists only of black vertices, 
at most one edge in each pair of crossing edges can be a chord of $C'$.
Therefore, we can take such a near-triangulation $T$ so that 
$C'$ has no chords (see the center of Fig.\ref{G'andH}). 
Observe that $T$ is $3$-connected. 
We denote by $F$ the unique non-triangular face of $T$. 

\begin{figure}[ht]
\centering
\input{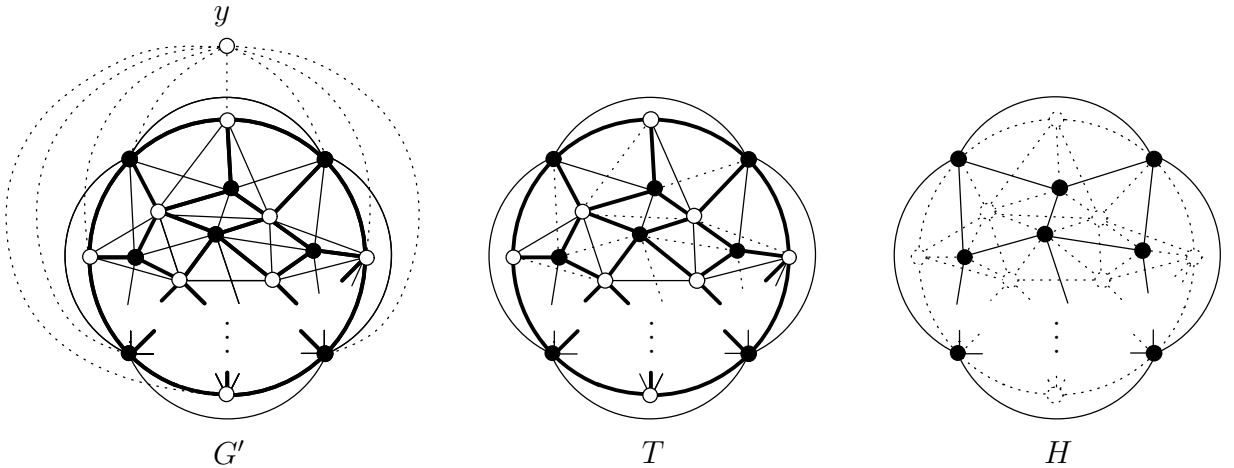}
\setlength{\abovecaptionskip}{10pt}
\caption{$G'$, $T$, and $H$}
\label{G'andH}
\end{figure}

Since $G$ is not $\{S_1, S_2, S_3\}$-linked,
there are no two disjoint paths 
joining $r_1$ to $r_2$ and $b_1$ to $b_2$ in $T$.
Hence, by Proposition \ref{seymour},
there exists an embedding of $T$ on $S^2$ 
such that there is a face whose boundary cycle 
containing $r_1, b_1, r_2$, and $b_2$ in this order. 
Since $T$ is $3$-connected, 
such an embedding is unique. 
Therefore, we may assume that 
$r_1, b_1, r_2,$ and $b_2$ are contained in $C'$ in this order. 

Let $H$ be the subgraph of $G$ induced by $B$ (see the right-hand side of Fig.\ref{G'andH}).
Observe that $H$ is a plane graph and the face $F$ of $T$ is 
also a face of $H$. 
Since $S_1, S_2 \subset V(C')$, we have
$S_1, S_2 \subset V(H)$. 
Furthermore, the boundary of each face of $H$ is a cycle, and hence $H$ is $2$-connected. 
Let $H_1$ (resp. $H_2$) 
be the graph obtained from $H$ by removing the vertices $r_1$ and $r_2$ (resp. $b_1$ and $b_2$). 
Suppose that $H_1$ or $H_2$, say $H_1$ here, is connected. 
Then, there exists a path $P_1$ joining $b_1$ and $b_2$ in $H_1$ such that $P_1$ 
only passes through the vertices of $B - \{r_1, r_2\}$.
Next, let $s$ (resp. $t$) be a vertex in $W \cap N_{G'}(r_1)$ (resp. $W \cap N_{G'}(r_2)$). 
Note that $s,t \notin \{b_1, b_2 \}$,
since $b_1, b_2 \in B$. 
Since the subgraph $G[W]$ of $G$ induced by $W$ is also $2$-connected, 
the subgraph $G[W]-y$ is connected. 
Then, there exists a path $P_{2}$ joining $s$ and $t$
such that $P_{2}$ only passes through the vertices of $W\setminus \{y\}$. 
Then, the path $P_{2} \cup \{r_1,r_2\} \cup \{r_1s, r_2t\}$
and $P_1$ are disjoint.
This contradicts the assumption. 
The same argument applies when $H_2$ is connected. 

Therefore, we assume that both $H_1$ and $H_2$ are disconnected. 
Since $H_1$ is disconnected, 
\{$r_1,r_2$\} is a $2$-cut of $H$. 
Then, there exists a separating $2$-curve $\gamma_1$
passing through only $r_1$ and $r_2$. 
Then $\gamma_1$ passes through the interior of a face $F_1$ of $H$ other than $F$ 
whose boundary contains $r_1$ and $r_2$. 
Since $H$ is the subgraph of $G$ induced by $B$,  
there is exactly one vertex $w\in W$ 
inside the face $F_1$. 
Similarly, there is a face $F_2$ of $H$ whose boundary 
contains $b_1$ and $b_2$ and whose interior contains exactly
one vertex $u \in W$. 
Under the situation, this implies that $F_1 = F_2$ whose boundary 
contains $r_1, b_1, r_2, b_2$ in this order, and that $w = u$. 
Therefore, $\{r_1,r_2,b_1,b_2,y\}$ have the arrangement in $G$ 
as stated in the Theorem \ref{thm:2-exh}. 
\end{Proofof}

Now, we can obtain the following corollary from Theorem \ref{thm:2-exh}.
\begin{cor}\label{(2,2,1)}
Every O1PG with connectivity $6$ is $(2, 2, 1)$-linked.
\end{cor}

\section{Concluding remarks and related topics}\label{sect:remarks}
\begin{figure}[ht]
\centering
\input{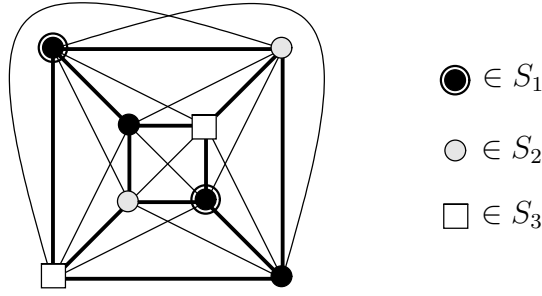}
\caption{$X\!W_6$ that is not $(2,2,2)$-linked}
\label{XW6}
\end{figure}

In the previous sections, we established the $(5,5)$-linkedness and the 
$(2,2,1)$-linkedness of O1PGs with connectivity $6$.
A natural question is whether similar results hold for other tuples of integers. 
In fact, there exist infinitely many O1PGs with connectivity 
$6$ that are not $(2,2,2)$-linked.
The following O1PG with eight vertices, shown in Fig.\ref{XW6}, called $X\!W_6$,
is one such example; 
there are no three vertex-disjoint connected subgraphs, 
each containing $S_i$ in the figure. 
As mentioned above, 
an $XW_6$ is a graph in an infinite series of graphs denoted by $X\!W_{2k}$ 
with $k\geq 3$ (see, e.g., \cite{k7minorssuzuki} for a definition). 
In fact, we have already observed that none of them is $(2,2,2)$-linked; 
since the required case analysis is extensive, a detailed discussion will be given in a subsequent paper. 

We briefly discuss several topics related to linkages.
Let $G$ be a graph and let $S$ be a vertex subset of $G$. 
The pair $(G, S)$ is {\em knitted} if for every partition of $S$ into non-empty subsets 
$S_1, S_2, \ldots, S_t$, there are disjoint connected subgraphs $G_1, G_2, \ldots, G_t$ 
in $G$ such that $S_i \subseteq V(G_i)$ for each $i \in \{1, \ldots,t \}$. 
A graph $G$ is $l$-{\em knitted} if 
$(G, S)$ is knitted for all $S \subseteq V(G)$ with $|S|=l$. 
The notion of a ``knitted graph'' was introduced by Bollob\'{a}s and Thomason \cite{bollo},
and was later further generalized by Kawarabayashi and Yu \cite{knitted}. 
Observe that a graph $G$ is $l$-knitted if and only if 
$G$ is $(s_1, s_2, \ldots, s_k)$-linked
for any list $(s_1,s_2, \ldots, s_k)$ 
with $\sum_{i=1}^k s_i = l$. 
In \cite{knitted}, the notion of a knitted graph was applied to study
the connectivity of minimum counterexamples to Hadwiger's conjecture.

By Corollary \ref{(2,2,1)} and brief arguments, 
we can obtain the following proposition. 

\begin{pro}
Every O1PG with connectivity $6$ is $5$-knitted. 
\end{pro}

However, there are O1PGs with connectivity $6$ that are not $6$-knitted, 
since there are infinitely many O1PGs with connectivity $6$ that are not $(2,2,2)$-linked 
as mentioned above; 
on the other hand, recall that every O1PG is $(3, 3)$-linked. 

We introduce one more related topic.
A graph $G$ is $k$-{\em ordered} if for every $k$ vertices of a given order, there is a cycle 
containing the $k$ vertices of the given order. 
Note that if a graph $G$ is $k$-ordered, then $G$ is $k$-knitted. 
Goddard \cite{Goddard4ordered} showed the following theorem. 

\begin{thm} [Goddard \cite{Goddard4ordered}] \label{4-ordered} 
Every $4$-connected triangulation on the sphere is $4$-ordered. 
\end{thm}

In fact, the result above has been extended to general closed surfaces 
by Mukae and Ozeki \cite{4ordered2010}. 
By Theorem \ref{4-ordered} and the fact that every O1PG contains a $4$-connected triangulation as a 
spanning subgraph, we obtain the following proposition. 

\begin{pro}
Every O1PG is $4$-ordered. 
\end{pro}

\begin{figure}[htbp]
\centering
\input{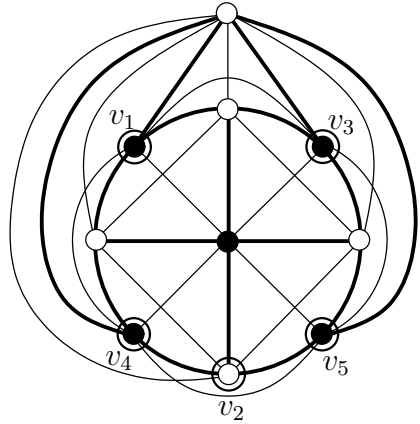}
\caption{$X\!W_8$ that is not $5$-ordered}
\label{XW8}
\end{figure}

However, there exist O1PGs with connectivity $6$ that are not $5$-ordered. 
The following O1PG called $X\!W_8$ is one such example (see Fig.\ref{XW8}). 
We can see that there are no cycles containing $v_1, v_2, v_3, v_4,$ and 
$v_5$ in this order. 

From the discussion so far, we can see that there is 
a clear gap among linkage, knittedness, and orderedness. 
Finally, we propose several open problems.

\begin{problem}
For a fixed integer $k$ and for any list $(s_1, s_2, \ldots, s_k)$ of positive integers, 
characterize $(s_1, s_2, \ldots, s_k)$-linked O1PGs. 
\end{problem}

\begin{problem}
Characterize $6$-knitted O1PGs with connectivity $6$. 
\end{problem}

\begin{problem}
Characterize $5$-ordered O1PGs with connectivity $6$. 
\end{problem}





\begin{thebibliography}{99}

\bibitem{bollo}
    B.~Bollob\'as and A.~Thomason,
    Highly linked graphs,
    {\em Combinatorica}
    {\bf 16}
    (1996), 313--320.    


\bibitem{chen}
    G.~Chen,
    R.~J.~Gould,
    K.~Kawarabayashi,
    F.~Pfender and B.~Wei,
    Graph minors and linkages,
    {\em Journal of Graph Theory}
    {\bf 49}
    (2005), 75--91.    

\bibitem{mnlinkedenami}
    K.~Enami and S.~Maezawa,
    Characterization of $(m, n)$-linked planar graphs,
    {\em Graphs and Combinatorics}
    {\bf 38}
    (2022), 131.
    
\bibitem{4n-8}
    I.~Fabrici and T.~Madaras,
    The structure of $1$-planar graphs,
    {\em Discrete Mathematics}
    {\bf 307}
    (2007), 854--865.
    
\bibitem{ConnectivitySuzuki}
    J.~Fujisawa,
    K.~Segawa and Y.~Suzuki,
    The matching extendability of 
    optimal $1$-planar graphs,
    {\em Graphs and Combinatorics}
    {\bf 34}
    (2018), 1089--1099.

\bibitem{Goddard4ordered}
   W.~Goddard,
   4-connected maximal planar graphs are 4-ordered,
   {\em Discrete Mathematics}
    {\bf 257}
    (2002), 405--410


\bibitem{klinked1}
    K.~Kawarabayashi,
    $k$-linked graphs with girth condition,
    {\em J. Graph Theory}
    {\bf 45}
    (2004), 48--50. 

\bibitem{knitted}
    K.~Kawarabayashi and G.~Yu
    Connectivities for {$k$}-knitted graphs and for minimal 
    counterexamples to {H}adwiger's conjecture,
    {\em J. Combin. Theory Ser. B}
    {\bf 103}
    (2013), 320--326.


\bibitem{Koizumi2025}
    S.~Koizumi and Y.~Suzuki,
    The matching extendability of
    optimal $1$-embedded graphs on the projective plane,
    {\em Discussiones Mathematicae Graph Theory}
    {\bf 45}
    (2025), 1233--1248.    

\bibitem{larman}
    D. G.~Larman and P.~Mani, 
    On the existence of certain configurations within graphs and 
    the {$1$}-skeletons of polytopes,
    {\em Proc. London Math. Soc. (3)}
    {\bf 20}
    (1970), 144--160.

\bibitem{mori200833linked}
    R.~Mori,
    $(3, 3)$-linked planar graphs,
    {\em Discrete Mathematics}
    {\bf 308}
    (2008), 5280--5283.

\bibitem{4ordered2010}
   R.~Mukae and K.~Ozeki,
   4-connected triangulations and 4-orderedness,
   {\em Discrete Mathematics}
    {\bf 310}
    (2010), 2271--2272

\bibitem{NNS} T.~Nagasawa, K.~Noguchi, Y.~Suzuki, 
Optimal $1$-embedded graphs which triangulate other surfaces, 
{\em J. Nonlinear Convex Anal.} {\bf 19} (2018), 1759--1770. 

\bibitem{relationship2015suzuki}
    K.~Noguchi and Y.~Suzuki,
    Relationship among triangulations, quadrangulations
    and optimal $1$-planar graphs,
    {\em Graphs and Combinatorics}
    {\bf 31}
    (2015), 1965--1972.

\bibitem{Ringel1965o1pg}
    G.~Ringel, 
    Ein Sechsfarbenproblem auf der kugel,
    {\em Abhandlungen aus dem Mathematischen Seminar der Universit\"at Hamburg}
    {\bf 29}
    (1965), 107--117.












\bibitem{Seymour1980}
    P.D.~Seymour,
    Disjoint paths in graphs,
    {\em Discrete Mathematics}
    {\bf 29}
    (1980), 293--309.

\bibitem{shiloach}
    Y.~Shiloach,
    A polynomial solution to the undirected two paths problem,
    {\em Journal of the ACM}
    {\bf 27}
    (1980), 445--456.


\bibitem{Re}
Y.~Suzuki, Re-embeddings of maximum $1$-planar graphs, 
{\em SIAM J. Discrete Math.} {\bf 24} (2010), 1527--1540.

\bibitem{k7minorssuzuki}
    Y.~Suzuki,
    $K_7$-minors in optimal $1$-planar graphs,
    {\em Discrete Mathematics}
    {\bf 340}
    (2017), 1227--1234. 
    
\bibitem{thomassen}
    C.~Thomassen,
    $2$-linked graphs,
    {\em European Journal of Combinatorics}
    {\bf 1}
    (1980), 371--378.







\end{thebibliography}
\end{document}